\documentclass{amsart}

\theoremstyle{definition}

\theoremstyle{remark}

\numberwithin{equation}{section}

\newcommand{\abs}[1]{\lvert#1\rvert}

\reversemarginpar

\newcommand{\nn}{\nonumber}
\newcommand{\no}{\noindent}

\newcommand{\realpart}{\mathop{\rm Re}\nolimits}

\newcommand{\lp}{\ln \sqrt{2 \pi}}
\newcommand{\lgf}{\ln \Gamma(\tfrac{1}{4})}
\newcommand{\ba}{\begin{eqnarray}}
\newcommand{\ea}{\end{eqnarray}}
\newcommand{\fq}{\frac{\partial}{\partial q}}
\newcommand{\ft}{\frac{\partial}{\partial t}}
\newcommand{\fnu}{\frac{\partial}{\partial \nu}}
\newcommand{\ift}{\int_{0}^{\infty}}
\newcommand{\ione}{\int_{0}^{1}}
\newcommand{\fmt}{\lfloor{\tfrac{m-1}{2}\rfloor}}
\newcommand{\fm}{\lfloor{\tfrac{m}{2}\rfloor}}

\newcommand{\allR}{\mathbb{R}}
\newcommand{\allN}{\mathbb{N}}
\newcommand{\nnN}{\mathbb{N}_{0}}
\newcommand{\st}{}

\newcommand{\strr}{ \left\{ k \atop j \right\} }
\newcommand{\stkj}{ \left\{ 2k \atop j \right\} }
\newcommand{\stkjo}{ \left\{ 2k+1 \atop j \right\} }

\begin{document}

\title[On some families of integrals] {On some families of integrals solvable in terms
of polygamma and negapolygamma functions}

\author{George Boros}
\address{Department of Mathematics,
Xavier University, New Orleans, LA 70125}
\email{gboros@xula.edu}

\author{Olivier Espinosa}
\address{Departamento de F\'{\i}sica,
Universidad T\'{e}cnica Federico Santa Mar\'{\i}a, Valpara\'{\i}so, Chile}
\email{espinosa@fis.utfsm.cl}

\author{Victor H. Moll}
\address{Department of Mathematics,
Tulane University, New Orleans, LA 70118}
\email{vhm@math.tulane.edu}

\subjclass{Primary 33}

\date{\today}

\keywords{Hurwitz zeta function, polygamma functions, loggamma, integrals}

\begin{abstract}
Beginning with Hermite's integral representation of the Hurwitz zeta
function, we derive explicit expressions in terms of elementary, polygamma, 
and negapolygamma functions for several families of integrals of
the type $\ift f(t)K(q,t)dt$ with kernels $K(q,t)$ equal to
$\left(e^{2\pi q t}-1\right)^{-1}$,
$\left(e^{2\pi q t}+1\right)^{-1}$, and
$\left(\sinh(2\pi q t)\right)^{-1}$.
\end{abstract}

\maketitle


\newtheorem{Definition}{\bf Definition}[section]
\newtheorem{Thm}[Definition]{\bf Theorem}
\newtheorem{Lem}[Definition]{\bf Lemma}
\newtheorem{Cor}[Definition]{\bf Corollary}
\newtheorem{Prop}[Definition]{\bf Proposition}
\newtheorem{Example}[Definition]{\bf Example}

\section{Introduction}
\label{S:intro}

The Hurwitz zeta function, defined by
\ba
\zeta(z,q) & = & \sum_{n=0}^{\infty} \frac{1}{(n+q)^{z}}
\ea
\no
for $z \in \mathbb{C}, \, \realpart{z} > 1,$ and $q \neq 0, -1, -2, \cdots,$
admits the
integral representation
\ba
\zeta(z,q) & = & \frac{1}{\Gamma(z)} \int_{0}^{\infty}
\frac{e^{-qt}}{1-e^{-t}} t^{z-1} dt,
\ea
\no
where $\Gamma(z)$ is Euler's gamma function, which is valid for
$\realpart z > 1$ and $\realpart q >0$, and can be used to
prove that $\zeta(z,q)$
admits an analytic extension to the whole complex plane except for a
simple pole at $z=1$. Hermite
proved an alternate integral representation, which actually
provides an explicit realization of this analytic continuation
for real $q>0$:
\ba
&  & \label{hermite1} \\
\zeta(z,q) & = & \frac{1}{2}q^{-z} + \frac{1}{z-1}q^{1-z} +
2q^{1-z} \int_{0}^{\infty} \frac{\sin(z \tan^{-1}t)}
{(1+t^{2})^{z/2} \, (e^{2 \pi t q} -1)}\, dt. \nn
\ea

\medskip

Special cases of $\zeta(z,q)$ include the  Bernoulli polynomials,
\ba
B_{m}(q) & = & -m \, \zeta(1-m,q),  \quad m \in \allN, \label{zetaber}
\ea
\no
defined by their generating function
\ba
\frac{te^{qt}}{e^{t}-1} & = & \sum_{m=0}^{\infty} B_{m}(q)
\frac{t^{m}}{m!}
\label{generating1}
\ea
\no
and given explicitly in terms of the Bernoulli numbers $B_{k}$ by
\ba
B_{m}(q) & = & \sum_{k=0}^{m} \binom{m}{k} B_{k}q^{m-k}, \label{berpoly}
\ea
and the polygamma functions,
\ba\label{hurwitz-polygamma}
\psi^{(m)}(q)= (-1)^{m+1}\,{m!}\,\zeta(m+1,q),\quad m\in\allN,
\ea
defined by
\ba \psi^{(m)}(q) & := & \frac{d^{\,m+1}}{dq^{\,m+1}} \ln\Gamma(q),\quad
m\in\allN.
\label{polygamma-def}
\ea
\medskip

The function $\zeta(z,q)$ is analytic for $z \neq 1$, and direct
differentiation of (\ref{hermite1}) yields
\begin{multline}\label{hermite2a}
\zeta'(z,q)  =  -\frac{1}{2}q^{-z} \ln q - \frac{q^{1-z}}{(z-1)^{2}}
- \frac{q^{1-z}}{z-1} \ln q
 -  2q^{1-z} \ln q \ift \frac{\sin(z \tan^{-1} t) \, dt}
{(1+t^{2})^{z/2} (e^{2 \pi q t} -1) }
\\
+  2q^{1-z} \ift \frac{\cos(z \tan^{-1} t)  \, \tan^{-1} t \; dt}
{(1+t^{2})^{z/2} (e^{2 \pi q t} -1) }
 -  q^{1-z} \ift \frac{\sin(z \tan^{-1} t)  \, \ln(1+t^{2}) \; dt}
{(1+t^{2})^{z/2} (e^{2 \pi q t} -1) },
\end{multline}
where $\zeta'(z,q)$ denotes $\partial\zeta(z,q)/\partial z$.
\medskip

In this paper we derive, starting from the representations
(\ref{hermite1}) and (\ref{hermite2a}), several definite integral
evaluations of the type
\[
F(q)=\ift \frac{f(t)}{e^{2 \pi q t} -1}\, dt.
\]
\no
The main examples considered here are the families
\ba
I_{k}(q) & = & \ift \frac{t}{(1+t^{2})^{k+1} (e^{2 \pi q t} -1) }\, dt, \nn \\
T_{k}(q) & = & \ift \frac{t^{k} \tan^{-1}t}{e^{2 \pi q t} -1 }\, dt, \nn \\
L_{k}(q) & = & \ift \frac{t^{k} \ln(1+t^{2})}{e^{2 \pi q t} -1 }\, dt, \nn
\ea
\no
and the associated integrals obtained by replacing the factor
$e^{2 \pi q t} -1$ in the denominator of the integrands by
$e^{2 \pi q t} + 1$ and $\sinh (2 \pi q t)$.
We produce closed-form expressions for $I_{k}(q)$ in terms of the polygamma
functions $\psi^{(m)}(q), \, 1\le m \le k$, and for $T_{2k}(q)$ and
$L_{2k+1}(q)$,
in terms of the derivative of the Hurwitz zeta function at
negative integers or, equivalently, the {\em balanced} functions
\ba\label{Ak-def}
A_{m}(q) & := & m\,\zeta'(1-m,q),
\ea
or the balanced negapolygamma functions,
\\
\ba\label{bal-negapolygamma}
\psi ^{( -m)} (q): = \frac{1}{m!}\left[{A_m (q)} - H_{m-1} B_m(q)\right],
\ea
\\
defined for $m\in\allN$, which were introduced in
\cite{esmo2}. $H_r$ is the harmonic number ($H_0:=0$).
We define a function $f(q)$ to be {\em balanced} (on the unit interval)
if it satisfies the properties
\[
\ione f(q)dq=0,\quad \text{and}\quad f(0)=f(1).
\]

For certain particular rational values of $q$ the balanced
negapolygamma functions evaluate to rational linear combinations
of elementary functions of special constants such as $\ln 2,
\ln\pi$, the Euler constant $\gamma$, $G/\pi$ ($G$ is Catalan's
constant), $\zeta'(-1)$, etc.
\\

We note that, in view of
Lerch's result \cite{hurber}
\ba
\ln \Gamma(q) & = & \zeta'(0,q) - \zeta'(0), \label{lerch}
\ea
\no
$A_{1}(q) = \zeta'(0,q)$ can be expressed in terms of the gamma function as
\ba
A_{1}(q) & = & \ln \frac{\Gamma(q)}{\sqrt{2 \pi} }.
\ea

\medskip

The problem of closed-form expressions for $T_{2k+1}$ and $L_{2k}$ remains
open. \\

\section{A series expansion}
\label{S:ramanujan}

All the results presented in this paper are consequences of the Taylor
series expansion of the function
\ba
f(z,t) & = & \frac{\sin(z \tan^{-1} t) }{(1+t^{2})^{z/2}},
\ea
\no
which appears in the integral representation of the Hurwitz zeta
function:
\ba
\zeta(z,q) &   = & \frac{1}{2}q^{-z} +
\frac{1}{z-1}q^{1-z} + 2q^{1-z} \int_{0}^{\infty}
\frac{f(z,t)}{e^{2 \pi tq}-1}\, dt.
\ea
\no

\begin{Thm}
The Taylor series
\ba
\frac{\sin(z \tan^{-1} t) }{(1+t^{2})^{z/2}} & = &
\sum_{k=0}^{\infty} \frac{(-1)^{k} (z)_{2k+1} }{(2k+1)!} t^{2k+1}
\label{expan1}
\ea
\no
and
\ba
\frac{\cos(z \tan^{-1} t) }{(1+t^{2})^{z/2}} & = &
\sum_{k=0}^{\infty} \frac{(-1)^{k} (z)_{2k} }{(2k)!} t^{2k}
\label{expan2}
\ea
\no
hold for $\abs{t}<1$.
\end{Thm}
\begin{proof}
Both sides of (\ref{expan1}) satisfy the equation
\ba
(1+t^{2}) \frac{d^{2}g}{dt^{2}} + 2t(z+1) \frac{dg}{dt} + z(z+1) g
= 0
\ea
\no
and the initial conditions $g(0)=0, \, g'(0) = z.$
It is straightforward to show that the series on
the right-hand side of (\ref{expan1}) converges for $\abs{t}<1$.
The proof of (\ref{expan2}) is similar.
\end{proof}

\medskip

\begin{Cor}
Let $m \in \mathbb{N}$. Then, for $t\in\allR$,
\ba
\cos(m \tan^{-1}t ) & = & (1+t^{2})^{-m/2}
\sum_{k=0}^{\fm} (-1)^{k} \binom{m}{2k} t^{2k} \label{cosm}
\ea
\no
and
\ba
\sin(m \tan^{-1}t ) & = & (1+t^{2})^{-m/2}
\sum_{k=0}^{\fmt} (-1)^{k} \binom{m}{2k+1} t^{2k+1}. \label{sinm}
\ea
\end{Cor}
\begin{proof}
The expression
\[
( - m)_n  = ( - 1)^n n!\binom{m}{n}  \st
\]
vanishes for $n>m$, so the series (\ref{expan1}, \ref{expan2})
terminate for $z=-m$.
\end{proof}

Since one can write $t^2=(1+t^2) - 1$, it is clear that, for
$m\in\allN$, the functions
$t^{-1} (1+t^2)^{m/2} \sin(m \tan^{-1} t)$ and
$(1+t^2)^{m/2} \cos(m \tan^{-1} t)$ are also polynomials in $1+t^{2}$.
We now give their explicit forms. \\

\begin{Cor}\label{sin-cos-expansion}
Let $m \in \mathbb{N}$. Then
\begin{multline}
\label{cosine2}
\cos(m \tan^{-1} t) = \sum_{p=0}^{\fm}
(-1)^{p} \frac{m}{m-p} \binom{m-p}{p} 2^{m-2p-1}
(1+t^{2})^{p-m/2} \st
\end{multline}
\no
and
\begin{multline}
\label{sine2}
\sin(m \tan^{-1} t) = t \, \sum_{p=0}^{\fmt}
(-1)^{p} \binom{m-p-1}{p} 2^{m-2p-1} (1+t^{2})^{p-m/2}. \st
\end{multline}
\no
\end{Cor}
\begin{proof}
Performing the binomial expansion of $t^{2k}=\left[(1+t^2) - 1\right]^k$
in (\ref{sinm}) we have
\begin{align}
\sin(m \tan^{-1}t) & = (1+t^{2})^{-m/2}
\sum_{k=0}^{\fmt}(-1)^{k} \binom{m}{2k+1} t^{2k+1} \st \nn \\
& = t(1+t^2)^{-m/2} \sum_{j=0}^{\fmt}
(-1)^{j} \left\{ \sum_{k=j}^{\fmt} \binom{m}{2k+1} \binom{k}{j}
\right\} (1+t^2)^{j}. \st \nn
\end{align}
The result now follows from the identity
\ba
\sum_{k=j}^{\fmt} \binom{m}{2k+1} \binom{k}{j} & = &
\binom{m-j-1}{j} 2^{m-2j-1}, \st
\ea
\no
where $0 \leq j \leq \fmt$.
A similar argument gives the identity for cosine.
\end{proof}

\bigskip

\section{A family of integrals derived from Hermite's representation}
\label{S:hermite}

The Hermite representation (\ref{hermite1}) can be written as
\ba
\int_0^\infty  {\frac{{\sin (z\tan ^{ - 1} t)}}
{{(1 + t^2 )^{z/2} (e^{2\pi qt}  - 1)}}} \,dt = \frac{1}
{2}q^{z - 1} \zeta (z,q) - \frac{1}
{{4q}} - \frac{1}
{{2(z - 1)}}.
\label{hermite1a}
\ea

\no
A direct consequence of the expansion (\ref{sinm}), when used in
(\ref{hermite1a}) with $z\in -\allN$, is the following 
well-known relation (\ref{zetaber})
between the Bernoulli polynomials and the Hurwitz zeta function.

\begin{Lem}
The Bernoulli polynomials $B_{m}(q)$, $m \in\allN$, are given by
\ba
B_{m}(q) & = & -m \,\zeta(1-m,q).
\ea
\end{Lem}
\begin{proof}
Substitute (\ref{sinm})
into (\ref{hermite1a}) with $z=-m$ and
use the well-known result \cite{gr} (3.411.2)
\ba
\int_0^\infty  {\frac{{t^{2k + 1} }}
{{e^{2\pi qt}  - 1}}} dt & = & ( - 1)^k \frac{{B_{2k + 2} }}
{{4(k + 1)q^{2k + 2} }}
\label{gr-34112}
\ea
\no
to obtain
\ba\label{bernoulli1}
\zeta ( - m,q) &=&  - \frac{{q^{m + 1} }}
{{m + 1}} + \frac{{q^m }}
{2} - 2q^{m + 1} \sum\limits_{k = 0}^{\left\lfloor {\frac{{m - 1}}
{2}} \right\rfloor } {\binom{m}
{{2k + 1}}\frac{{B_{2k + 2} }}
{{4(k + 1)q^{2k + 2} }}}\nn\\
&=& - \frac{1}
{{m + 1}}\sum\limits_{k = 0}^{m + 1} {\binom{{m + 1}}
{k}B_k q^{m + 1 - k} } \nn \\
&=&- \frac{1}{{m + 1}}B_{m+1}(q), \nn
\ea
\no
in view of (\ref{berpoly}) and the facts that 
$B_0=1, B_1=-1/2,$ and $B_{2k+1}=0$ for $k\in\allN$.
\end{proof}

\medskip

Similarly, the alternate expansion (\ref{sine2}), 
when substituted into
(\ref{hermite1a}) with $z=m+1,\,m\in\allN$,
leads us to
consider the following family of integrals. 

\begin{Thm}
The integrals
\ba
\label{in-def0}
I_{k}(q) & := & \ift \frac{t}{(1+t^2)^{k+1} (e^{2 \pi q t} -1) }\, dt,
\ea
$k\in\allN$, are given  by
\begin{multline}
\label{In-explicit0}
I_{k} (q) = -\frac{1}{4k} - \frac{\binom{2k}{k} }{2^{2k+2}q}
+\frac{1}{k 2^{2k}} \sum_{j=1}^{k}
\frac{(-1)^{j+1}}{(j-1)!} \binom{2k-j-1}{k-j} 2^{j-1} q^{j}
\psi^{(j)}(q). \st
\end{multline}
\end{Thm}
\begin{proof}
Use the expansion (\ref{sine2}) in (\ref{hermite1a}) with $z=m+1$ to obtain
the recursion
\ba
\sum\limits_{p = 0}^{\left\lfloor {\frac{{m}}
{2}} \right\rfloor } ( - 1)^p 2^{m-2p+1}  \binom{m-p}{p}
I_{m - p} (q) = q^{m} \zeta (m+1,q) - \frac{1}
{{2q}} - \frac{1}{m} \nn
\ea
\no
for $m\in\allN$. Expressing the Hurwitz zeta functions in terms of
polygamma functions by inverting (\ref{hurwitz-polygamma}), we find
\begin{multline}\label{recursion1}
\sum\limits_{p = \left\lfloor {\frac{{m + 1}}{2}}
\right\rfloor }^m {( - 1)^p 2^{2p} \binom{p}
{{m - p}}} I_{p} (q) =  - 2^{m-1} \left[ {\frac{{q^m \psi ^{(m)} (q)}}
{{m!}} + \frac{{( - 1)^m }}
{{2q}} + \frac{{( - 1)^m }}
{m}} \right]. \st
\end{multline}
\no
The recursion (\ref{recursion1}) can be solved in closed
form by inverting the sum. We first lower the lower limit to 
$p=1$ since the binomial coefficient vanishes when
$p < m-p$, and then use the orthogonality formula

\ba
\sum\limits_{j = 1}^k ( - 1)^j j \binom{2k - j - 1}{k-j} \binom{p}{j-p}
=
\begin{cases}
(-1)^k k  & \text{ if } p=k  \\
0 & \text{otherwise}
\end{cases} \st
\nn
\ea
and the evaluations
\ba
\sum_{j=1}^{k} 2^{j} \binom{2k-j-1}{k-j}  = 2^{2k-1} & \text{ and } &
\sum_{j=1}^{k} j \, 2^{j} \binom{2k-j-1}{k-j}  = k \binom{2k}{k}
\nn
\ea
to obtain the explicit formula  (\ref{In-explicit0}).

\end{proof}

\medskip

\no
{\bf Note}. The case $k=0$ appears in \cite{gr} (3.415.1):
\ba
\label{I1-explicit}
I_0 (q) = \int_0^\infty  {\frac{t}
{{(1+t^2)(e^{2\pi qt}  - 1)}}}\, dt = \frac{1}{2} \ln q - \frac{1}
{{4q}} - \frac{1}{2} \psi (q),
\ea
where $\psi (q)$ is the digamma function.
This result also follows from (\ref{hermite1a}) in the limit $z\to 1$,
in view of
\begin{align}
\psi(q)&=\lim_{z\to1}\left[ \frac{1}{z-1}-\zeta(z,q)\right].
\label{psi-hurwitz}
\end{align}
\bigskip

\section{Two new families  of integrals}
\label{S:2families}

As we know from Corollary \ref{sin-cos-expansion},
the functions
\[
t^{-1} (1+t^2)^{-z/2} \, \sin(z \tan^{-1} t)
\quad\text{and}\quad
(1+t^2)^{-z/2} \cos(z \tan^{-1} t)
\]
\noindent
are polynomials in $1+t^{2}$
when $z \in - \allN$, a fact that allows considerable
simplification of representation (\ref{hermite2a}). This leads us
to consider the families of integrals
\ba
T_{k}(q) &  = & \ift \frac{t^{k} \tan^{-1} t}{e^{2 \pi q t} -1} \, dt
\label{T}
\ea
\no
and
\ba
L_{k}(q) &  = & \ift \frac{t^{k} \ln(1+t^{2})}{e^{2 \pi q t} -1} \,
dt,
\ea
\no
that appear after differentiating
Hermite's representation (\ref{hermite1}) with respect to the
parameter $z$. A direct differentiation of (\ref{hermite1a})
yields
\begin{multline}\label{hermite2}
\int_0^\infty  {\frac{{\cos (z\tan ^{ - 1} t)\tan ^{ - 1} t}}
{{(1 + t^2 )^{z/2} (e^{2\pi qt}  - 1)}}} \,dt - \frac{1}
{2}\int_0^\infty  {\frac{{\sin (z\tan ^{ - 1} t)\ln (1 + t^2 )}}
{{(1 + t^2 )^{z/2} (e^{2\pi qt}  - 1)}}} \,dt \\= \frac{1}
{2}\left[ {q^{z - 1} \zeta (z,q)\ln q + q^{z - 1} \zeta' (z,q)
+ \frac{1}{{(z - 1)^2 }}} \right].
\end{multline}
Setting $z=-m, \,m\in\nnN$, we find a recursion for these integrals. \\

\begin{Prop}
For $m\in\nnN$, the integrals $T_{2k}(q)$ and $L_{2k+1}(q)$ satisfy the relation
\begin{multline}\label{JplusK}
2 \sum_{k=0}^{\fm} (-1)^{k} \binom{m}{2k}T_{2k}(q) +
\sum_{k=0}^{\fmt} (-1)^{k} \binom{m}{2k+1}L_{2k+1}(q)
\\
 = \frac{1}{(m+1)q^{m+1}}\left[ A_{m+1}(q) - B_{m+1}(q)\,\ln q\right] +
\frac{1}{(m+1)^{2}}.
\end{multline}
\end{Prop}
\begin{proof}
This is derived directly from (\ref{hermite2})
using the expansions (\ref{cosm}) and (\ref{sinm}), and the relations
(\ref{zetaber}) and (\ref{Ak-def}).
\end{proof}

\medskip

Equation (\ref{JplusK}) can be used iteratively to find explicit
expressions for the integrals $T_{2k}(q)$ and $L_{2k+1}(q)$ in terms of
the functions $A_m(q)$ and $B_m(q)$.  \\

\begin{Example}
\label{ex1}
The value $m=0$ in (\ref{JplusK}) yields
\[
T_0 (q) = \frac{1} {2q}A_1 (q) - \frac{1} {2q}B_1 (q)\ln q + \frac{1}{2}.
\]
Using (\ref{lerch}) and $B_{1}(q) = q - \tfrac{1}{2}$, we have
\begin{equation}
\label{T0:explicit}
T_{0}(q) = \ift \frac{\tan^{-1}t}{e^{2 \pi q t} -1}\, dt  =
\frac{1}{2} - \frac{1}{2}\ln q + \frac{\ln q} {4q} + \frac{{\ln \Gamma
(q)}} {2q} - \frac{{\lp}} {2q}.
\end{equation}
We note that this result corresponds to Binet's second expression
for $\ln\Gamma(q)$ \cite{ww}.

\medskip

\no Some particular evaluations of $T_{0}(q)$ are
\[
\begin{gathered}
  \int_0^\infty  {\frac{{\tan ^{ - 1} t}}
{{e^{2\pi t}  - 1}}\, dt}  = \frac{1}{2}  - \frac{\lp}{2}, \hfill \\
  \int_0^\infty  {\frac{{\tan ^{ - 1} t}}
{{e^{\pi t}  - 1}}\, dt}  = \frac{1}{2} - \frac{\ln 2}{2}, \hfill \\
  \int_0^\infty  {\frac{{\tan ^{ - 1} t}}
{{e^{\pi t/2}  - 1}}\, dt}  = \frac{1}{2} - \ln \pi -2 \ln 2 + 2 \lgf.  \hfill \\
\end{gathered}
\]
\end{Example}

\medskip

\begin{Example}
The case $m=1$ in (\ref{JplusK}) yields
\ba
2T_{0}(q) + L_{1}(q) & = & \frac{A_{2}(q)}{2q^{2}} -
\frac{B_{2}(q) \, \ln q}{2q^{2}} + \frac{1}{4},
\ea
\no
and using $A_2(q)=2\zeta'(-1,q)$, $B_{2}(q) = q^2 - q + 1/6,$ and 
the expression for $T_{0}(q)$
obtained in (\ref{T0:explicit}),  we get

\ba
L_1 (q) & = & \ift \frac{t \, \ln (1+t^{2})}{e^{2 \pi qt} -1} \, dt  \label{l1}
\\
& =  & \frac{1} {{q^2 }}\zeta'(-1,q) - \frac{{\ln \Gamma (q)}} {q} +
\frac{{\lp}} {q} - \left( {\frac{1} {{12q^2 }} - \frac{1} {2}}
\right)\ln q - \frac{3} {4}. \nn
\ea
\no

We see that $L_{1}(q)$ will evaluate to special values whenever
$\zeta'(-1,q)$ does. Particular examples of the latter are
\begin{align}
\zeta'(-1,\tfrac{1}{2})&=-\tfrac{1}{2}\zeta '( -
1)-\tfrac{1}{24}\ln2\nn\\
\intertext{and}
\zeta'(-1,\tfrac{1}{4})&=-\tfrac{1}{8}\zeta '( - 1)+ G/4 \pi,
\nn
\end{align}
\\
\no
and particular evaluations of
$L_{1}(q)$ are
\[
\begin{gathered}
  \int_0^\infty  {\frac{{t\ln (1 + t^2 )}}
{{e^{2\pi t}  - 1}}\, dt}  = \zeta '( - 1) + \lp - \frac{3}
{4}, \hfill \\
  \int_0^\infty  {\frac{{t\ln (1 + t^2 )}}
{{e^{\pi t}  - 1}}\, dt}  =  - 2\zeta '( - 1) + \frac{2}
{3}\ln 2 - \frac{3}
{4}, \hfill \\
  \int_0^\infty  {\frac{{t\ln (1 + t^2 )}}
{{e^{\pi t/2}  - 1}}\, dt}  =  - 2\zeta '( - 1) + \frac{5}
{3}\ln 2 - \frac{3}
{4} + \frac{{4G}}
{\pi } - 4 \lgf
+ 4\lp. \hfill \\
\end{gathered}
\]
\end{Example}

\medskip

\bigskip

We now evaluate the integrals $T_{2k}(q)$ and $L_{2k+1}(q)$
in terms of elementary functions and the
balanced negapolygamma functions (\ref{bal-negapolygamma}).  \\

\begin{Thm}
\label{tandl}
Let $k \in \allN$. Then,
\ba
(-1)^{k} T_{2k}(q) & = & (-1)^{k} \ift \frac{t^{2k} \tan^{-1} t}{e^{2 \pi q t} -1} \, dt \label{t2k} \\
 & = &  \frac{1}{2(2k+1)^{2}} - \frac{\ln q}{2(2k+1)}
+ \frac{1}{8kq}  \nn \\
& & +
\frac{1}{4} \sum_{j=0}^{k-1} \frac{B_{2j+2}}{(j+1)(2k-2j-1)}  \frac{1}{q^{2j+2}} \nn \\
 & & + \frac{1}{2} \sum_{j=0}^{2k} (-1)^j \frac{(2k)!}{(2k-j)!}
\frac{\psi^{(-1-j)}(q)}{q^{j+1}}   \nn
\ea
\no
and
\ba
T_0 (q) & = & 
\frac{1}{2} - \frac{\ln q}{2} + \frac{\ln q} {4q} + \frac{{\ln \Gamma
(q)}} {2q} - \frac{{\lp}} {2q}. \nn
\ea
\no
Similarly, for $k \geq 0$,
\begin{align}
\label{l2kp1}
(-1)^{k+1}L_{2k+1}(q) & =  (-1)^{k+1} \ift \frac{t^{2k+1} \ln (1+t^2) }
{e^{2 \pi q t } -1 } \, dt  \\
 & = \frac{1}{(2k+2)^{2}} - \frac{\ln q}{2k+2} +
\frac{1}{2q(2k+1)}  \nn \\
&\quad + \frac{1}{2} \sum_{j=0}^{k-1} \frac{B_{2j+2}}{(j+1)(2k-2j)} \frac{1}{
q^{2j+2}} + \frac{B_{2k+2}}{(2k+2)q^{2k+2} } \left( \ln q - H_{2k+1} \right)
\nn \\
 &\quad + \sum_{j=0}^{2k+1} (-1)^j \frac{(2k+1)!}{(2k-j+1)!}
\frac{\psi^{(-1-j)}(q)}{q^{j+1}}.  \nn
\end{align}
\end{Thm}
\begin{proof}
We first derive the expression for $T_{2k}(q)$. For any
differentiable function $f$  we have
\ba
q \frac{\partial}{\partial q} f(qt) & = &
t \frac{\partial}{\partial t} f(qt), \nn
\ea
so
\ba
q \frac{d}{dq} T_{2k}(q) & = & \ift \left( t^{2k} \tan^{-1}t \right)
q \fq \left( \frac{1}{e^{2 \pi q t } -1} \right) dt \nn \\
 & = & \ift \left( t^{2k+1} \tan^{-1} t \right) \ft
\left( \frac{1}{e^{2 \pi q t } -1} \right) dt \nn \\
& = & -(2k+1)T_{2k}(q) - \ift \frac{t^{2k+1} \; dt }
{(e^{2 \pi q t } -1)(1+t^{2})}. \nn
\ea
\no
Thus
\ba
\frac{d}{dq} \left( q^{2k+1} T_{2k}(q) \right) & = &
-q^{2k} \ift \frac{t^{2k+1} \; dt}{(e^{2 \pi q t} -1) (1+t^{2})},  \nn
\ea
\no
and since
\ba
t^{2k} & = & (-1)^{k} + (-1)^{k+1} (1+t^{2}) \; \sum_{j=0}^{k-1}
(-1)^{j} t^{2j}, \nn
\ea
\no
we have
\begin{multline}
(-1)^{k+1} \frac{d}{dq} \left( q^{2k+1} T_{2k}(q) \right) =
q^{2k} \ift \frac{t \, dt}{(e^{2 \pi q t } -1)(1+t^{2})}
\\
- q^{2k} \sum_{j=0}^{k-1} (-1)^{j} \ift \frac{t^{2j+1} \, dt}
{(e^{2 \pi q t}-1)}. \nn
\end{multline}
\no
Using (\ref{gr-34112}) and (\ref{I1-explicit})  we obtain
\[
(-1)^{k+1} \frac{d}{dq} \left( q^{2k+1} T_{2k}(q) \right)  =
\frac{1}{2} q^{2k} \ln q - \frac{1}{4}q^{2k-1} - \frac{1}{2}q^{2k} \psi(q)
- \frac{1}{4} \sum_{j=0}^{k-1} \frac{B_{2j+2}}{j+1} q^{2k-2j-2}.
\]

Now, for $k \ge 1$ each of the terms on the right-hand side is integrable
at $q=0$, so that
\ba
& & \label{T2n-explicit1}  \\
(-1)^{k+1} T_{2k}(q) & = &  \frac{-1}{2(2k+1)^{2}} + \frac{\ln q}{2(2k+1)}
- \frac{1}{8kq} - \frac{1}{2q^{2k+1}} \int_{0}^{q} r^{2k} \psi(r) dr \nn \\
 & & - \frac{1}{4} \sum_{j=0}^{k-1} \frac{B_{2j+2}}{j+1}
\frac{q^{-2j-2}}{2k-2j-1} + \frac{c_{2k}}{q^{2k+1}}, \nn
\ea
\no
where $c_{2k}$ is a constant of integration which can be determined by
studying the behavior of $q^{2k+1} T_{2k}(q)$ as $q \to 0$:

\ba
c_{2k}   & =  & ( - 1)^{k+1}  \lim _{q \to 0} q^{2k + 1} T_{2k} (q) \nn \\
& = &( - 1)^{k+1}  \frac{{(2k)!}}
{{4(2\pi )^{2k} }}\zeta (2k + 1) = -\frac{1}{2} \zeta '( - 2k).
\nn
\ea
The evaluation of the limit above is obtained by replacing $\tan^{-1}(x/q)$ by
$\pi/2$ in
\[
q^{2k + 1} T_{2k} (q) = \int_0^\infty  {\frac{{x^{2k} \tan ^{ - 1} (x/q)}}
{{e^{2\pi x}  - 1}}dx}
\]
and employing formula \cite{gr}(3.411.1):
\ba\label{gr-34111}
\int_0^\infty  {\frac{{x^{\nu  - 1} }}
{{e^{\mu x}  - 1}}dx}  = \frac{1}
{{\mu ^\nu  }}\Gamma (\nu )\zeta (\nu ),
\quad\realpart\mu  > 0,\realpart\nu  > 1.
\ea

The final step in the evaluation of $T_{2k}(q)$ uses the result 
\begin{multline}
\label{intpoly1}
\int_{0}^{q} r^{n} \psi(r)\,dr = n!
\sum_{j=0}^{n} \frac{(-1)^j}{(n-j)!} q^{n-j}\psi^{(-1-j)}(q)
-n! (-1)^n \psi^{(-1-n)}(0),
\end{multline}
\no
valid for $n\in\allN$, which can be obtained from the corresponding
indefinite integral given in \cite{esmo2}.
Since
\[
\psi^{(-1-n)}(0) = \frac{1}{n!}\left[
\zeta'(-n)-\frac{H_{n}B_{n+1}}{n+1}\right],
\]
we see that for $n=2k$ the boundary term above precisely cancels
the term proportional to the integration constant $c_{2k}$, thus
leading to the explicit formula (\ref{t2k}).

\medskip

The formula for $L_{2k+1}(q)$ is derived in a similar way.  We start with
\ba
\ln (1+t^{2}) & = &
\frac{d}{dt} \left[ t \ln (1+t^{2}) -2t + 2 \tan^{-1}t \right] \nn
\ea
\no
and integrate by parts, observing that
\ba
\ft \left[ \frac{t^{2k+1}}{e^{2 \pi q t} -1 } \right] & = &
\frac{(2k+1)t^{2k}}{e^{2 \pi q t } -1} + t^{2k} q \fq
\left[ \frac{1}{e^{2 \pi q t} -1 } \right],  \nn
\ea
\no
to conclude that
\ba
\fq \left( q^{2k+2} \, L_{2k+1}(q) \right) & = & -2q \fq \left( q^{2k+1}
T_{2k}(q) \right) +
\frac{(-1)^{k+1} B_{2k+2} }{(2k+2)q }. \nn
\ea
\no
Using the expression (\ref{T2n-explicit1}) for $T_{2k}(q)$ we obtain
$L_{2k+1}(q)$ up to a constant of integration:
\begin{multline}
(-1)^{k+1} L_{2k+1}(q) = \frac{1}{(2k+2)^{2}} - \frac{\ln q}{2k+2}
+ \frac{1}{2q(2k+1)} + \frac{1}{q^{2k+2}} \int_{0}^{q} r^{2k+1} \psi(r)dr
\\
+ \frac{1}{2} \sum_{j=0}^{k-1} \frac{B_{2j+2} }{(j+1)(2k-2j)}
\frac{1}{q^{2j+2}}  + \frac{B_{2k+2} \, \ln q}{(2k+2)q^{2k+2}}  +
\frac{c_{2k+1}}{q^{2k+2}}.
\nn
\end{multline}
\no
As before, the constant $c_{2k+1}$ can be determined by evaluating
\ba
c_{2k+1} &  = & \lim\limits_{q \to 0} \left[(-1)^{k+1} q^{2k+2} L_{2k+1}(q) -
\frac{B_{2k+2}}{2k+2} \ln q\right]. \nn
\ea
\no
Note that
\ba
q^{2k+2} L_{2k+1}(q) & = & -2 \ln q \ift \frac{x^{2k+1} \; dx}{e^{2 \pi x} -1 }
+ \ift \frac{x^{2k+1} \, \ln(x^{2} + q^{2}) }{e^{2 \pi x} -1 } dx,
\nn
\ea
\no
so that, in view of (\ref{gr-34112}), the limit above is given simply by
\ba
(-1)^{k+1} c_{2k+1}  & = &  2 \ift \frac{x^{2k+1} \, \ln x}{e^{2 \pi x}-1}dx
\nn\\
& = & 2 \fnu \left[ \frac{ \Gamma( \nu) \zeta(\nu) }{(2 \pi)^{\nu}}
\right] \Bigg{|}_{\nu = 2k+2}
=  \fnu \left[ \frac{\zeta(1 - \nu) }{\cos( \pi \nu/2) }
\right]\Bigg{|}_{\nu=2k+2}  \nn \\
& = & (-1)^{k} \zeta'(-2k-1), \nn
\ea
\no
and thus $c_{2k+1} = - \zeta'(-2k-1)$. In the second line above
we used the functional equation for the Riemann zeta function.
This time, however,  the boundary term from (\ref{intpoly1}) and the term
containing the integration constant $c_{2k+1}$ cancel only partially,
leaving the term proportional to the harmonic number $H_{2k+1}$ that
appears in (\ref{l2kp1}).\\
\end{proof}

\no
{\bf Note}. Adamchik \cite{adam2}
informed us that he is able to evaluate the same
integrals in terms of the Barnes G-function.  This function is uniquely
defined by the recurrence formula
\ba
G_{n+1}(z+1) & = & \frac{G_{n+1}(z)}{G_{n}(z)}, \label{barnes1} \\
G_{1}(z) & = & 1/\Gamma(z), \nn
\ea
\no
and the condition
\ba
 \frac{d^{n+1}}{dx^{n+1}} \left\{ \log G_{n}(x) \right\} & \geq  & 0.
\ea
\no
The integrals $T_{2k}(q)$ and $L_{2k+1}(q)$ are reportedly given by
\begin{multline}
(-1)^k T_{2k}(q) =
-\frac{\ln q - H_{2k+1}}{2(2k+1)} + \frac{1}{2}  \sum_{j=0}^{2k} (-1)^{j}
\binom{2k}{j} \zeta'(-j)q^{-j-1}\\
- \frac{1}{2 q^{2k+1}}
 \sum_{j=1}^{2k} (-1)^{j} j! \stkj \ln G_{j+1}(q+1)
\nn
\end{multline}
\no
and
\begin{multline}
(-1)^{k+1} L_{2k+1}(q) = \frac{B_{2k+2} }{2k+2} \, \frac{\ln q}{q^{2k+2}}
 - \frac{\ln q - H_{2k+2}}{2k+2} \\
+\sum_{j=0}^{2k+1} (-1)^{j} \binom{2k+1}{j} \zeta'(-j) q^{-j-1}\\
- \frac{1}{q^{2k+2}} \sum_{j=1}^{2k+1} (-1)^{j} j!
\stkjo \ln G_{j+1}(q+1),
\nn
\end{multline}
\no
where $\strr$ are the Stirling numbers of the second kind and $H_{k}$ are the
harmonic numbers.  \\

We have been unable to
determine the values of $T_{2k+1}(q)$ and $L_{2k}(q)$ using the techniques
described here.

\medskip

\section{Some related integrals}
\label{additional}

The formulas for definite integrals developed in the previous sections
involve the kernel $(e^{2 \pi q t}-1)^{-1}$. These evaluations, combined
with a simple manipulation, lead to a larger class. \\

\begin{Lem}
Let
\ba
F(q) & = & \ift \frac{f(t)}{e^{2 \pi q t } -1}\, dt.
\ea
\no
Then
\ba
G(q) & := & \ift \frac{f(t)}{e^{2 \pi q t} +1}\, dt =  F(q) - 2 F(2q),
\\
S(q) & : = & \ift \frac{f(t)}{\sinh(2 \pi qt) }\, dt = 2F(q) - 2F(2q).
\ea
\end{Lem}
\begin{proof}
This is a direct consequence of the identities
\ba
\frac{1}{e^{x}+1} & = & \frac{1}{e^{x}-1} - \frac{2}{e^{2x}-1}, \nn \\
\frac{1}{\sinh x } & = & \frac{2}{e^{x}-1} - \frac{2}{e^{2x}-1}. \nn
\ea
\end{proof}

\no
\begin{Example}
The expression (\ref{T0:explicit}) yields
\ba
\ift \frac{\tan^{-1}t}{e^{2 \pi q t}+1} \, dt & = &
 \ln 2 - \frac{1}{2} + \frac{\ln q}{2} - \frac{\ln 2}{4q}
+ \frac{\ln \Gamma(q) - \ln \Gamma(2q) }{2q},
\ea
\no
with special evaluations
\ba
\ift \frac{\tan^{-1} t}{e^{2 \pi t } + 1}\, dt  & = &
\frac{3}{4} \ln 2 - \frac{1}{2}, \nn \\
\ift \frac{\tan^{-1} t}{e^{ \pi t } + 1}\, dt  & = &
\frac{\ln \pi }{2} - \frac{1}{2}, \nn \\
\ift \frac{\tan^{-1} t}{e^{ \pi t/2 } + 1}\, dt  & = &
- \frac{1}{2}  -\ln 2 + 2 \lgf - \ln \pi. \nn
\ea
\no
Similarly
\ba
\ift \frac{\tan^{-1} t}{\sinh(2\pi q t)}\, dt & = &
\ln 2 - \frac{\ln(4 \pi)}{4q} + \frac{\ln q}{4q}
+ \frac{\ln \Gamma(q)}{q} - \frac{\ln \Gamma(2q)}{2q}.
\ea

\medskip

\no
Some particular values are
\[
\begin{gathered}
\ift \frac{\tan^{-1} t}{\sinh(2 \pi t)}\, dt  =
\frac{1}{2} \ln 2 - \frac{1}{4} \ln \pi, \hfill \\
\ift \frac{\tan^{-1} t}{\sinh(\pi t)}\, dt  =
\frac{1}{2} \ln \pi - \frac{1}{2} \ln2, \hfill \\
\ift \frac{\tan^{-1} t}{\sinh(\pi t/2)}\, dt =
4 \lgf  - 2 \ln \pi - 3 \ln 2. \hfill \\
\end{gathered}
\]
\end{Example}

\no
\begin{Example}
The expression (\ref{In-explicit0}) yields
\begin{multline}\label{plusin}
\ift \frac{t}{(1+t^{2})^{k+1} \, (e^{2 \pi q t} +1) }\, dt  =
\frac{1}{4k} \\
+  \frac{1}{k2^{2k+1}}
\sum_{j=1}^{k} \frac{(-1)^{j+1}}{(j-1)!} \binom{2k-j-1}{k-j} 2^{j}q^{j}
\left[ \psi^{(j)}(q) - 2^{j+1} \psi^{(j)}(2q) \right]
\end{multline}
\no
and
\begin{multline}\label{sinhin}
\ift \frac{t}{(1+t^{2})^{k+1} \, \sinh(2 \pi q t)}\, dt =
- \frac{1}{2^{2k+2}q} \binom{2k}{k} \\
+ \frac{1}{k2^{2k}}
\sum_{j=1}^{k} \frac{(-1)^{j+1}}{(j-1)!} \binom{2k-j-1}{k-j} 2^{j}q^{j}
\left[ \psi^{(j)}(q) - 2^{j} \psi^{(j)}(2q) \right].
\end{multline}
\end{Example}

\medskip

\no
The function $\psi(q)$ and its derivatives do not satisfy a simple duplication
formula. Thus the explicit evaluation of (\ref{plusin}) and (\ref{sinhin})
requires the values of $\psi^{(j)}$ at $q$ and $2q$.

\medskip

\begin{Example}
The expression (\ref{t2k}) yields
\ba
& & \label{mess1} \\
(-1)^k \ift \frac{t^{2k} \, \tan^{-1}t}{e^{2 \pi q t}+1}\, dt &  =  &
\frac{-1}{2(2k+1)^2} + \frac{\ln 2}{2k+1} + \frac{\ln q}{2(2k+1)} \nn  \\
& + & \frac{1}{4} \sum_{j=0}^{k-1}
\frac{B_{2j+2} (1 - 2^{-2j-1})}{(j+1)(2k-2j-1)q^{2j+2}} \nn\\
 & + & \frac{1}{2} \sum_{j=0}^{2k} \frac{(-1)^{j} (2k)!}{(2k-j)!q^{j+1}}
\left[ \psi^{(-1-j)}(q) - \frac{\psi^{(-1-j)}(2q)}{2^{j}} \right] \nn
\ea
\no
and
\ba
& & \label{mess2} \\
(-1)^k \ift \frac{t^{2k} \tan^{-1} t}{\sinh (2 q \pi t) }\, dt & = &
\frac{\ln 2}{2k+1} + \frac{1}{8kq} \nn \\
& + & \frac{1}{2} \sum_{j=0}^{k-1}
\frac{B_{2j+2} (1 - 2^{-2j-2})}{(j+1)(2k-2j-1)q^{2j+2}} \nn\\
 & + &  \sum_{j=0}^{2k} \frac{(-1)^{j} (2k)!}{(2k-j)!q^{j+1}}
\left[ \psi^{(-1-j)}(q) - \frac{\psi^{(-1-j)}(2q)}{2^{j+1}} \right]. \nn
\ea

\end{Example}

\begin{Example}
The expression (\ref{l2kp1}) yields
\begin{multline}
(-1)^{k+1} \ift \frac{t^{2k+1} \, \ln(1+t^2)}{e^{2 \pi q t} +1}\, dt  =
- \frac{1}{(2k+2)^{2}} + \frac{2\ln 2}{2k+2} + \frac{\ln q}{2k+2}
\\
 +  \frac{1}{2} \sum_{j=0}^{k-1}
\frac{B_{2j+2} \, (1 - 2^{-2j-1} )}{(j+1)(2k-2j) q^{2j+2}}
\\
+  \frac{B_{2k+2}}{(2k+2)q^{2k+2} }
\left[ ( 1 - \frac{1}{2^{2k+1}} ) \ln q - \frac{ \ln 2 }{2^{2k+1}} -
(1 - \frac{1}{2^{2k+1}} ) H_{2k+1} \right]
\\
 +   \sum_{j=0}^{2k+1} \frac{(-1)^{j} (2k+1)! }{(2k-j+1)!q^{j+1}}
\left[ \psi^{(-1-j)}(q) - \frac{\psi^{(-1-j)}(2q)}{2^{j}}
\right]
\end{multline}
\no
and
\begin{multline}
(-1)^{k+1} \ift \frac{t^{2k+1} \, \ln(1+t^2)}
{\sinh (2 \pi q t) }\, dt = \frac{2\ln 2}{2k+2} + \frac{1}{2q(2k+1)}
\\
+ \sum_{j=0}^{k-1} \frac{B_{2j+2} (1 - 2^{-2j-2})}{(j+1)(2k-2j)q^{2j+2}}
\\
+ \frac{2B_{2k+2}}{(2k+2)q^{2k+2}}
\left[ (1 - \frac{1}{2^{2k+2}}) \ln q - \frac{\ln 2}{2^{2k+2}}
- (1 - \frac{1}{2^{2k+2}}) H_{2k+1} \right]
\\
 +  2  \sum_{j=0}^{2k+1} \frac{(-1)^{j} (2k+1)! }{(2k-j+1)!q^{j+1}}
\left[ \psi^{(-1-j)}(q) - \frac{\psi^{(-1-j)}(2q)}{2^{j+1}}\right].
\end{multline}

\no
For example,
\begin{multline}
\ift \frac{t \, \ln(1+t^2)}{\sinh(2 \pi q t)}\, dt =
\frac{1}{q^{2}} \left[ 2\zeta'(-1, q) - \frac{1}{2}\zeta'(-1, 2q) \right]
- \, \frac{1}{q} \left( 2 \ln \Gamma(q) - \ln \Gamma(2 q) \right) \\
 - \ln 2 + \frac{\ln \pi}{2q} - \frac{\ln q}{8q^{2}}
 + \frac{\ln 2}{24q^{2}} + \frac{\ln 2}{2q}.
\end{multline}

\medskip

\no
Some special values are

\ba
\ift \frac{t \, \ln(1+t^2)}{\sinh(2 \pi t) } \, dt & = &
-\frac{11}{24} \ln 2 + \frac{1}{2} \ln \pi + \frac{3}{2} \zeta'(-1),  \nn \\
\ift \frac{t \, \ln(1+t^2)}{\sinh( \pi t) } \, dt & = &
\frac{1}{3} \ln 2 - \ln \pi - 6 \zeta'(-1), \nn \\
\ift \frac{t \, \ln(1+t^2)}{\sinh( \pi t/2) } \, dt & = &
6 \ln 2 + 4\ln \pi + \frac{8G}{\pi} - 8 \lgf. \nn
\ea
\end{Example}

\no
{\bf Note}. 
Differentiating with respect to the parameter $q$ and evaluating at special
values yields many new integrals.  For example, the derivative of
(\ref{T0:explicit}) at $q=1$ and $q=2$ yields
\ba
\ift \frac{t \, \tan^{-1} t}{\sinh^{2} \pi t} \, dt & = & \frac{1}{2 \pi}
+ \frac{\gamma}{\pi}  - \frac{\ln \sqrt{2 \pi}}{\pi}
\ea
\no
and
\ba
\ift \frac{t \, \tan^{-1} t}{\sinh^{2} 2 \pi t} \, dt & = & -\frac{1}{8 \pi}
+ \frac{\gamma}{2 \pi}  - \frac{\ln \pi}{8 \pi},
\ea
\no
respectively. 

\medskip

\no
{\bf Acknowledgments}. The second author would like to thank the Department of
Mathematics at Tulane University for its hospitality. The
third author acknowledges the
partial support of NSF-DMS 0070567, Project number 540623.

\end{document}